\documentclass[10pt,a4paper]{article}
\usepackage{latexsym,fullpage,amsfonts,amssymb,amsmath,amscd,graphics,epic,theorem}
\usepackage[all]{xy}
\usepackage{amsmath}
\usepackage{amssymb, amstext}
\usepackage{hyperref}

\baselineskip 18pt \textwidth 16cm  \theoremstyle{plain}
\newtheorem{lemma}{Lemma}[subsection]
\newtheorem{proposition}[lemma]{Proposition}
\newtheorem{remark}[lemma]{Remark}
\newtheorem{example}[lemma]{Example}
\newtheorem{theorem}[lemma]{Theorem}
\newtheorem{definition}[lemma]{Definition}
\newtheorem{notation}[lemma]{Notation}

\newtheorem{property}{}
\newtheorem{corollary}[lemma]{Corollary}

\begin{document}
\newcommand{\pperp}{\hbox{$\perp\hskip-6pt\perp$}}
\newcommand{\N}{{\mathbb N}}
\newcommand{\PP}{{\mathbb P}}
\newcommand{\Z}{{\mathbb Z}}
\newcommand{\Q}{{\mathbb Q}}
\newcommand{\R}{{\mathbb R}}
\newcommand{\C}{{\mathbb C}}
\newcommand{\K}{{\mathbb K}}
\newcommand{\F}{{\mathbb F}}
\newcommand{\proofend}{\hfill$\Box$\smallskip}
\newcommand{\eps}{{\varepsilon}}
\newcommand{\ko}{{\mathcal O}}
\newcommand{\wx}{{\widetilde x}}
\newcommand{\wz}{{\widetilde z}}
\newcommand{\wa}{{\widetilde a}}
\newcommand{\bz}{{\boldsymbol z}}
\newcommand{\bp}{{\boldsymbol p}}
\newcommand{\wy}{{\widetilde y}}
\newcommand{\wc}{{\widetilde c}}
\newcommand{\bi}{{\omega}}
\newcommand{\bx}{{\boldsymbol x}}
\newcommand{\Log}{{\operatorname{Log}}}
\newcommand{\pr}{{\operatorname{pr}}}
\newcommand{\Graph}{{\operatorname{Graph}}}
\newcommand{\jet}{{\operatorname{jet}}}
\newcommand{\Tor}{{\operatorname{Tor}}}
\newcommand{\sqh}{{\operatorname{sqh}}}
\newcommand{\const}{{\operatorname{const}}}
\newcommand{\Arc}{{\operatorname{Arc}}}
\newcommand{\Sing}{{\operatorname{Sing}}}
\newcommand{\Span}{{\operatorname{Span}}}
\newcommand{\Aut}{{\operatorname{Aut}}}
\newcommand{\Ker}{{\operatorname{Ker}}}
\newcommand{\Int}{{\operatorname{Int}}}
\newcommand{\Aff}{{\operatorname{Aff}}}
\newcommand{\Area}{{\operatorname{Area}}}
\newcommand{\val}{{\operatorname{Val}}}
\newcommand{\conv}{{\operatorname{conv}}}
\newcommand{\rk}{{\operatorname{rk}}}
\newcommand{\ow}{{\overline w}}
\newcommand{\ov}{{\overline v}}
\newcommand{\ks}{{\cal S}}
\newcommand{\red}{{\operatorname{red}}}
\newcommand{\kc}{{\cal C}}
\newcommand{\ki}{{\cal I}}
\newcommand{\kj}{{\cal J}}
\newcommand{\ke}{{\cal E}}
\newcommand{\kz}{{\cal Z}}
\newcommand{\tet}{{\theta}}
\newcommand{\Del}{{\Delta}}
\newcommand{\bet}{{\beta}}
\newcommand{\mm}{{\mathfrak m}}
\newcommand{\kap}{{\kappa}}
\newcommand{\del}{{\delta}}
\newcommand{\sig}{{\sigma}}
\newcommand{\alp}{{\alpha}}
\newcommand{\Sig}{{\Sigma}}
\newcommand{\Gam}{{\Gamma}}
\newcommand{\gam}{{\gamma}}
\newcommand{\Lam}{{\Lambda}}
\newcommand{\lam}{{\lambda}}
\newcommand{\om}{{\omega}}
\newcommand{\udot}{\mathaccent\cdot\cup}
\newcommand{\Sc}{{\cal S}}
\newcommand{\G}{{\cal G}}
\newcommand{\T}{{\cal T}}
\newcommand{\Fre}{{Fr\'{e}chet \,}}
\newcommand{\Che}{{Ch\v{e}ch \,}}
\newcommand{\ctp}{{\widehat{\otimes}}}

\title{De-Rham theorem and Shapiro lemma for Schwartz functions on Nash manifolds}
\author{Avraham Aizenbud and Dmitry Gourevitch \thanks{ Avraham Aizenbud and Dmitry Gourevitch, Faculty of Mathematics and Computer
Science, The Weizmann Institute of Science POB 26, Rehovot 76100,
ISRAEL. 
E-mails: aizenr@yahoo.com, guredim@yahoo.com.\smallskip
}   }
\date{}

\maketitle

\begin{abstract}
In this paper we continue our work on Schwartz functions and
generalized Schwartz functions on Nash (i.e. smooth
semi-algebraic) manifolds. Our first goal is to prove analogs of
de-Rham theorem for de-Rham complexes with coefficients in
Schwartz functions and generalized Schwartz functions. Using that
we compute cohomologies of the Lie algebra ${\mathfrak g}$ of an
algebraic group $G$ with coefficients in the space of generalized
Schwartz sections of $G$-equivariant bundle over a $G$- transitive
variety $M$. We do it under some assumptions on topological
properties of $G$ and $M$. This computation for the classical case
is known as Shapiro lemma.
\end{abstract}

\tableofcontents
\section{Introduction}

We will use the notions of Schwartz sections and generalized
Schwartz sections of Nash (i.e. smooth semi-algebraic) bundles
over Nash manifolds introduced in \cite{AG}. These will be
reviewed in section \ref{Prel}.

We use the following notations. For a Nash manifold $M$ we denote
by $\Sc(M)$ the space of Schwartz functions on $M$ and by $\G(M)$
the space of generalized Schwartz functions on $M$. For a Nash
vector bundle $E \rightarrow M$ we denote by $\Sc_M^E$ the cosheaf
of Schwartz sections of $E$ and by $\G_M^E$ the sheaf of
generalized Schwartz sections of $E$. We also denote the global
Schwartz sections of $E$ by $\Sc(M,E)$ and global generalized
Schwartz sections of $E$ by $\G(M,E)$.

Let $M$ be a Nash manifold. We can define the de-Rham complex with
coefficients in Schwartz functions $$DR_{\Sc}(M): 0 \rightarrow
{\Sc}(M,{\Omega^0_M}) \rightarrow ... \rightarrow
{\Sc}(M,{\Omega^n_M}) \rightarrow 0$$ We will prove that its
cohomologies are isomorphic to compact support cohomologies of
$M$. Similarly we will define de-Rham complex with coefficients in
generalized Schwartz functions $$DR_{\G}(M): 0 \rightarrow
\G(M,{\Omega^0_M}) \rightarrow ... \rightarrow \G(M,{\Omega^n_M})
\rightarrow 0$$ and prove that its cohomologies are isomorphic to
cohomologies of $M$.

Moreover, we will prove relative versions of these statements. Let
$F \rightarrow M$ be a locally trivial fibration. We will define
Nash vector bundles $H^i(F \rightarrow M)$ and $H^i_c(F
\rightarrow M )$ over $M$ such that their fibers will be equal to
the cohomologies of the fibers of $F$ and the compact support
cohomologies of the fibers of $F$ in correspondance. We will
define relative de-Rham complex of $F \rightarrow M$ with
coefficients in Schwartz functions. We will denote it by
$DR_{\Sc}(F \rightarrow M)$ and prove that its cohomologies are
canonically isomorphic to global Schwartz sections of the bundles
$H^i_c(F \rightarrow M)$.

Similarly we will define relative de-Rham complex  of $F
\rightarrow M$ with coefficients in generalized Schwartz functions
and denote it by $DR_{\G}(F \rightarrow M)$. We will prove that
its cohomologies are canonically isomorphic to global generalized
Schwartz sections of the bundles $H^i(F \rightarrow M )$. In
particular, if the fiber of $F$ is contractible then  the higher
cohomologies of the relative de-Rham complex with coefficients in
generalized Schwartz functions vanish and the zero cohomology is
$\G(M)$. Using this result we will prove the following analog of
Shapiro lemma.
\begin{theorem} 
Let $G$ be a contractible linear algebraic group. Let $H<G$ be a
contractible subgroup and let $M = G/H$. Let $\rho$ be a finite
dimensional representation of $H$. Let $E \rightarrow M$ be the
$G$-equivariant bundle corresponding to $\rho$. Let ${\mathfrak
h}$ be the Lie algebra of $H$ and ${\mathfrak g}$ be the Lie
algebra of $G$. Let $V$ be the space of generalized Schwartz
sections of $E$ over $M$. It carries a natural action of $G$.

Then the cohomologies of ${\mathfrak g}$ with coefficients in $V$
are isomorphic to the cohomologies of ${\mathfrak h}$ with
coefficients in $\rho$.
\end{theorem}
We will need Nash analogs of some known notions and theorems from
algebraic topology that we have not found in the literature. They
are written in section \ref{algtop}.

\subsection{Structure of the paper}

In section \ref{Prel} we give the necessary preliminaries on Nash
manifolds and Schwartz functions and distributions over them.

In subsection \ref{SemiAlg} we introduce basic notions of
semi-algebraic geometry from \cite{BCR}, and \cite{Shi}. In
particular we formulate the Tarski-Seidenberg principle of
quantifier elimination.

In subsection \ref{sheaftheo} we introduce the notion of
restricted topological space (from \cite{DK}) and sheaf theory
over it. These notions will be necessary to introduce non-affine
Nash manifolds and to formulate the relative de-Rham theorem.

In subsection \ref{Nash} we give basic preliminaries on Nash
manifolds from \cite{BCR}, \cite{Shi} and \cite{AG}.

In subsection  \ref{algtop} we repeat known notions and theorems
from algebraic topology for the Nash case. In particular we
formulate Theorem \ref{sursec} which says that the restricted
topology is equivalent as a Grothendieck topology to the smooth
topology on the category of Nash manifolds.

In subsection \ref{Schwartz} we give the definitions of Schwartz
functions and Schwartz distributions on Nash manifolds from
\cite{AG}.

In subsection \ref{NucFre} we remind some classical facts on
nuclear \Fre spaces and prove that the space of Schwartz functions
on a Nash manifold is nuclear.

In section \ref{DeRham} we formulate and prove de-Rham theorem for
Schwartz functions on Nash manifolds. Also, we prove its relative
version. We need this relative version in the proof of Shapiro
Lemma.

In section \ref{secShapLem} we formulate and prove a version of
Shapiro lemma for Schwartz functions on Nash manifolds.

In section \ref{Summary} we discuss possible extensions and
applications of this work.

In appendix \ref{AppSurSec} we prove Theorem \ref{sursec} that we
discussed above.

\subsection*{Acknowledgements} We would like to thank our
teacher {\bf Joseph Bernstein} for teaching us almost all of
mathematics we know and for his help in this work.

We would like to thank {\bf Semyon Alesker}, {\bf Lev Buhovski},
{\bf Vadim Kosoy}, {\bf  Vitali Milman}, {\bf Omer Offen}, {\bf
Lev Radzivilovski}, {\bf Eitan Sayag},  and {\bf  Eugenii Shustin}
for helpful discussions and {\bf Vladimir Berkovich}, {\bf Paul
Biran} and {\bf Sergei Yakovenko} for useful remarks.

\section{Preliminaries} \label{Prel}
During the whole paper we mean by smooth infinitely
differentiable.

\subsection{Semi-algebraic sets and Tarski-Seidenberg principle} \label{SemiAlg}

In this subsection we will give some preliminaries on
semi-algebraic geometry from \cite{BCR} and \cite{Shi}.

\begin{definition}
{  A subset $A \subset \R^n$ is called a \textbf{semi-algebraic
set} if it can be presented as a finite union of sets defined by a
finite number of polynomial equalities and inequalities. In other
words, if there exist finitely many polynomials $f_{ij}, g_{ik}
\in R[x_1,...,x_n]$ such that $$A = \bigcup \limits _{i=1}^r \{x
\in \R^n | f_{i1}(x)>0,...,f_{is_i}(x)>0,
g_{i1}(x)=0,...,g_{it_i}(x)=0\}.$$ }
\end{definition}
\begin{lemma}
The collection of semi-algebraic sets is closed with respect to
finite unions, finite intersections and complements.
\end{lemma}
\begin{example}
The semi-algebraic subsets of $\R$ are unions of finite number of
intervals.
\end{example}
\begin{proposition}
Let $\nu$ be a bijective semi-algebraic mapping. Then the inverse
mapping $\nu^{-1}$ is also semi-algebraic.
\end{proposition}
\emph{Proof.} The graph of $\nu$ is obtained from the graph of
$\nu^{-1}$ by switching the coordinates. $\text{ }$ \proofend

One of the main tools in the theory of semi-algebraic spaces is
the Tarski-Seidenberg principle of quantifier elimination. Here we
will formulate and use a special case of it. We start from the
geometric formulation.
\begin{theorem}
Let $A \subset \R^n$ be a semi-algebraic subset and $p:\R^n
\rightarrow \R^{n-1}$ be the standard projection. Then the image
$p(A)$ is a semi-algebraic subset of $\R^{n-1}$.
\end{theorem}
By induction and a standard graph argument we get the following
corollary.
\begin{corollary}
An image of a semi-algebraic subset of $\R^n$ under a
semi-algebraic map is semi-algebraic.
\end{corollary}
Sometimes it is more convenient to use the logical formulation of
the Tarski-Seidenberg principle. Informally it says that any set
that can be described in semi-algebraic language is
semi-algebraic. We will now give the logical formulation and
immediately after that define the logical notion used in it.
\begin{theorem} {(Tarski-Seidenberg principle)}
Let $\Phi$ be a formula of the language $L(\R)$ of ordered fields
with parameters in $\R$. Then there exists a quantifier - free
formula $\Psi$ of $L(\R)$ with the same free variables
$x_1,\dots,x_n$ as $\Phi$ such that
 $\forall x \in \R^n, \Phi(x) \Leftrightarrow \Psi(x)$.
\end{theorem}
For the proof see Proposition 2.2.4 on page 28 of \cite{BCR}.
\begin{definition}
{A \textbf{formula of the language of ordered fields with
parameters in $\R$} is a formula written with a finite number of
conjunctions, disjunctions, negations and universal and
existential quantifiers ($\forall$ and $\exists$) on variables,
starting from atomic formulas which are formulas of the kind
$f(x_1,\dots,x_n) = 0$ or $g(x_1,\dots,x_n) > 0$, where $f$ and
$g$ are polynomials with coefficients in $\R$. The free variables
of a formula are those variables of the polynomials which are not
quantified. We denote the language of such formulas by $L(\R)$.}
\end{definition}
\begin{notation}
Let $\Phi$ be a formula of $L(\R)$ with free variables
$x_1,\dots,x_n$. It defines the set of all points
$(x_1,\dots,x_n)$ in $\R^n$ that satisfy $\Phi$. We denote this
set by $S_{\Phi}$. In short, $$S_{\Phi}:=\{ x
\in \R^n |\Phi(x)\}.$$ 
\end{notation}
\begin{corollary}
Let $\Phi$ be a formula of $L(\R)$. Then $S_{\Phi}$ is a
semi-algebraic set.
\end{corollary}
\emph{Proof.} Let $\Psi$ be a quantifier-free formula equivalent
to $\Phi$. The set $S_{\Psi}$ is semi-algebraic since it is a
finite union of sets defined by  polynomial equalities and
inequalities. Hence $S_{\Phi}$ is also semi-algebraic since
$S_{\Phi}=S_{\Psi}$. \proofend
\begin{proposition}
The logical formulation of the Tarski-Seidenberg principle implies
the geometric one.
\end{proposition}
\emph{Proof.} Let $A \subset \R^n$ be a semi-algebraic subset, and
$pr:\R^n \rightarrow \R^{n-1}$ the standard projection. Then there
exists a formula $\Phi \in L(\R)$ such that $A = S_{\Phi}$. Then
$pr(A)= S_{\Psi}$ where
$$\Psi(y)= \text{``}\exists x \in \R^n \,(\pr(x) = y \wedge \Phi(x))".$$ Since
$\Psi \in L(\R)$, the proposition follows now from the previous
corollary.
\begin{remark}
In fact, it is not difficult to deduce the logical formulation
from the geometric one.
\end{remark}
Let us now demonstrate how to use the logical formulation of the
Tarski-Seidenberg principle.
\begin{corollary}
The closure of a semi-algebraic set is semi-algebraic.
\end{corollary}
\emph{Proof.} Let $A \subset \R^n$ be a semi-algebraic subset, and
let $\overline{A}$ be its closure. Then $\overline{A}=S_{\Psi}$
where
$$\Psi(x)= \text{``} \forall \eps>0 \, \exists y \in A \, |x-y|^2<\eps".$$
Clearly, $\Psi \in L(\R)$ and hence $\overline{A}$ is
semi-algebraic. \proofend
\begin{corollary}
Images and preimages of semi-algebraic sets under semi-algebraic
mappings are semi-algebraic.
\end{corollary}
%
\begin{corollary}
$ $\\(i) The composition of semi-algebraic mappings is
semi-algebraic.\\
(ii) The $\R$-valued semi-algebraic functions on a semi-algebraic
set $A$ form a ring, and any nowhere vanishing semi-algebraic
function is invertible in this ring.
\end{corollary}

We will also use the following theorem from \cite{BCR}
(Proposition 2.4.5).
\begin{theorem} \label{fincomp}
Any semi-algebraic set in $\R^n$ has a finite number of connected
components.
\end{theorem}

\subsection{Sheaf theory over restricted topological spaces}
\label{sheaftheo}
The usual notion of topology does not fit semi-algebraic geometry.
Therefore we will need a different notion of topology called
restricted topology, that was introduced in \cite{DK}.
\begin{definition}
{ A \textbf{restricted topological space $M$} is a set $M$
equipped with a family $ \overset {_{\circ}} {\mathfrak S}   (M)$
of subsets of $M$, called the open subsets that includes $M$ and
the empty set and is closed with respect to finite unions and
finite intersections.}
\end{definition}
\begin{remark}
In general, there is no closure in restricted topology  since
infinite intersection of closed sets does not have to be closed.
\end{remark}
\begin{remark}
A restricted topological space $M$ can be considered as a site in
the sense of Grothendieck. The category of the site has as objects
the open sets of $M$ and as morphisms the inclusion maps. The
covers $(U_i \rightarrow U)_{i\in I}$ are the \underline{finite}
systems of inclusions with $\bigcup \limits _{i=1}^n U_i = U$.
This gives us the notions of sheaf and cosheaf on $M$. We will
repeat the definitions of this notion in simpler terms.
\end{remark}
\begin{definition}
{  A \textbf{pre-sheaf $F$} on a restricted topological space $M$
is a contravarinant functor from the category $Top(M)$ which has
open sets as its objects and inclusions as morphisms to the
category of abelian groups, vector spaces etc.

In other words, it is the assignment $U \mapsto F(U)$ for every
open $U$ with abelian groups, vector spaces etc. as values, and
for every inclusion of open sets $V \subset U$ - a restriction
morphism $res_{U,V}: F(U) \rightarrow F(V)$ that satisfy
$res_{U,U} = Id$ and for $W \subset V \subset U$, $res_{V,W} \circ
res_{U,V} = res_{U,W}$. A \textbf{morphism of pre-sheaves} $\phi:
F \rightarrow G$ is a collection of morphisms $\phi_U:F(U)
\rightarrow G(U)$ for any open set $U$ that commute with the
restrictions.}
\end{definition}
%
\begin{definition}
{A \textbf{sheaf $F$} on a restricted topological space $M$  is a
pre-sheaf fulfilling the usual sheaf conditions, except that now
only finite open covers are admitted. The conditions are: for any
open set $U$ and any finite cover $U_i$ of $M$ by open subsets,
the sequence $$0 \rightarrow F(U) \overset{res_1}{\rightarrow}
\prod_{i=1}^n F(U_i) \overset{res_2}{\rightarrow}
\prod_{i=1}^{n-1} \prod_{j=i+1}^{n} F(U_i \cap U_j)$$ is exact.\\
The map $res_1$ above is defined by $res_1(\xi) = \prod \limits
_{i=1}^n res_{U,U_i}(\xi)$ and the map $res_2$ by
$$res_2(\prod_{i=1}^n \xi_i) = \prod_{i=1}^{n-1} \prod_{j=i+1}^n res_{U_i,U_i \cap U_j}(\xi_i)
- res_{U_j,U_i \cap U_j}(\xi_j)$$ . }
\end{definition}

\begin{definition}
{  A \textbf{pre-cosheaf $F$} on a restricted topological space
$M$ is a covarinant functor from the category $Top(M)$ to the
category of abelian groups, vector spaces etc.

In other words, it is the assignment $U \mapsto F(U)$ for every
open $U$ with abelian groups, vector spaces etc. as values, and
for every inclusion of open sets $V \subset U$ - an extension
morphism $ext_{V,U}: F(V) \rightarrow F(U)$ that satisfy:
$ext_{U,U} = Id$ and for $W \subset V \subset U$, $ext_{V,U} \circ
ext_{W,V} = ext_{W,U}$. A morphism of pre-cosheaves $\phi: F
\rightarrow G$ is a collection of morphisms $\phi_U:F(U)
\rightarrow G(U)$ for any open set $U$ that commute with the
extensions. }
\end{definition}


\begin{definition}
{A \textbf{cosheaf $F$} on a restricted topological space $M$ is a
pre-cosheaf on $M$ fulfilling the conditions dual to the usual
sheaf conditions, and with only finite open covers allowed. This
means: for any open set $U$ and any finite cover $U_i$ of $M$ by
open subsets, the sequence $$ \bigoplus_{i=1}^{n-1}
\bigoplus_{j=i+1}^{n} F(U_i \cap U_j) \rightarrow
\bigoplus_{i=1}^n F(U_i) \rightarrow F(U) \rightarrow 0$$ is
exact. \\
Here, the first map is defined by
$$ \bigoplus_{i=1}^{n-1} \bigoplus_{j=i+1}^{n} \xi_{ij} \mapsto \sum_{i=1}^{n-1} \sum_{j=i+1}^{n}
 ext_{U_i \cap U_j,U_i}(\xi_{ij}) -
 ext_{U_i \cap U_j,U_j}(\xi_{ij}) $$
and the second one by $$\bigoplus_{i=1}^n \xi_i \mapsto \sum
\limits  _{i=1}^n ext_{U_i,U}(\xi_i).$$  }
\end{definition}

\begin{remark}
As in the usual case, we have the functors of sheafification and
cosheafification, which assign to every pre-sheaf (pre-cosheaf) a
canonical sheaf (cosheaf). They are defined as left adjoint (right
adjoint) functors to the forgetful functor from sheaves
(cosheaves) to pre-sheaves (pre-cosheaves). Note that in the
construction of cosheafification quotient objects are needed. So
cosheafification always exists for sheaves with values in abelian
categories. Pre-cosheaves of \Fre spaces whose extension maps have
closed image also have cosheafification.
\end{remark}
\begin{definition}\label{suppdist}
{  Let $M$ be a restricted topological space, and $F$ be a sheaf
on $M$. Let $Z \subset M$ be a closed subset. A global section of
$F$ is said to be \textbf{supported in $Z$} if its restriction to
the complement of $Z$ is zero.}
\end{definition}
\begin{remark}
Unfortunately, if we will try to define support of a section, it
will not be a closed set in general, since infinite intersection
of closed sets in the restricted topology does not have to be
closed.
\end{remark}
\begin{remark}
Till the end of this section we will consider only sheaves and
cosheaves of linear spaces over $\R$.
\end{remark}

\begin{definition}
{  Let $M$ be a restricted topological space and $V$ be a linear
space over $\R$ . A function $f:M \rightarrow V$ is called
\textbf{locally constant} if there exists a \underline{finite}
cover $M = \bigcup \limits _{i=1}^k U_i$ s.t. $\forall i. f|_{U_i}
=const$.}
\end{definition}

\begin{remark}
Till the end of this section we will consider only those
restricted topological spaces in which any open set is a finite
disjoint union of its open connected subsets. In such spaces a
locally constant function is a function which is constant on every
connected component.
\end{remark}
Using this notion, we define constant sheaf in the usual way, i.e.

\begin{definition}
{  Let $M$ be a restricted topological space. Let $V$ be a linear
space over $\R$. We define \textbf{constant sheaf over $M$ with
coefficients in $V$} by $V_M(U):=\{f:U \rightarrow V| f$ is
locally constant on $V$ in the induced restricted topology $\}$
for any open $U\subset M$.}
\end{definition}

\begin{definition}
{  Let $M$ be a restricted topological space and $F$ be a sheaf
(cosheaf) over it. We define a conjugate cosheaf (sheaf) by
$F^*(U):=F(U)^*$. }
\end{definition}

\begin{definition}
{  Let $M$ be a restricted topological space. Let $V$ be a finite
dimensional linear space over $\R$. We define \textbf{constant
cosheaf over $M$ with coefficients in $V$} by $V_M':=({V^*_M})^*$
.}
\end{definition}

\begin{definition}
{  A sheaf(cosheaf) $F$ over a restricted topological space $M$ is
called \textbf{locally constant} if there exists a
\underline{finite} cover $M = \bigcup \limits _{i=1}^k U_i$ such
that for any $i$, $F|_{U_i}$ is isomorphic to a constant
sheaf(cosheaf) on $U_i$.}
\end{definition}

\begin{definition}
{  We define internal $\mathcal{H}om$ in the categories of sheaves
and cosheaves over restricted topological space the same way as it
is done in the usual case, i.e.
$\mathcal{H}om(F,G)(U):=Hom(F|_U,G|_U)$.}
\end{definition}

\begin{definition}
{  Let $F$ be a sheaf over a restricted topological space $M$. We
define its dual sheaf $D(F)$ by $D(F):=\mathcal{H}om(F,{\R}_M)$.}
\end{definition}

\begin{definition}
{  Let $F$ be a cosheaf over a restricted topological space $M$.
We define its dual cosheaf $D(F)$ by
$D(F):=\mathcal{H}om(F,\R_M')$.}
\end{definition}

\begin{notation}
{  To every sheaf(cosheaf) $F$ over a restricted topological space
$M$ we associate a cosheaf (sheaf) $F'$ by $F':=D(F)^*$ .}
\end{notation}

\begin{remark}
The constant sheaf (cosheaf) is evidently a sheaf (cosheaf) of
algebras, and any sheaf (cosheaf) has a canonical structure of a
sheaf (cosheaf) of modules over the constant sheaf (cosheaf).
\end{remark}

\begin{definition}
{  Let $F$ and $G$ be sheaves (cosheaves). We define $F\otimes G$
to be the sheafification (cosheafification) of the presheaf
(precosheaf) $U \mapsto F(U) \underset {\R_M(U)} {\otimes} G(U).$}
\end{definition}

\subsection{Nash manifolds} \label{Nash}
In this section we define the category of Nash manifolds,
following \cite{BCR}, \cite{Shi} and \cite{AG}.
\begin{definition}
{  A \textbf{Nash map} from an open semi-algebraic subset $U$ of
$\R^n$ to an open semi-algebraic subset $V\subset \R^m$ is a
smooth (i.e. infinitely differentiable) semi-algebraic function.
The ring of $\R$-valued Nash functions on $U$ is denoted by ${\cal
N}(U)$. A \textbf{Nash diffeomorphism} is a Nash bijection whose
inverse map is also Nash. }
\end{definition}
As we are going to do semi-algebraic differential geometry, we
will need a semi-algebraic version of implicit function theorem.
\begin{theorem}[Implicit Function Theorem.]
{  Let $(x^0,y^0) \in \R^{n+p}$, and let $f_1,...,f_p$ be
semi-algebraic smooth functions on an open neighborhood of
$(x^0,y^0)$, such that $f_j(x^0,y^0)=0$ for $j=1,..,p$ and the
matrix $[\frac{\partial f_j}{\partial y_i}(x^0,y^0)]$ is
invertible. Then there exist open semi-algebraic neighborhoods $U$
(resp. V) of $x^0$ (resp. $y^0$) in $\R^n$ (resp. $\R^p$) and a
Nash mapping $\phi$, such that $\phi(x^0)=y^0$ and $f_1(x,y) = ...
= f_p(x,y)=0 \Leftrightarrow y = \phi(x)$ for every $(x,y) \in U
\times V.$}
\end{theorem}
The proof is written on page 57 of \cite{BCR} (corollary 2.9.8).
\begin{definition}
{  A \textbf{Nash submanifold} of $\R^n$ is a semi-algebraic
subset of $\R^n$ which is a smooth submanifold .}
\end{definition}
By the implicit function theorem it is easy to see that this
definition is equivalent to the following one, given in
\cite{BCR}:
\begin{definition}
{  A semi-algebraic subset $M$ of $\R^n$ is said to be a
\textbf{Nash submanifold of $\R^n$ of dimension $d$} if, for every
point $x$ of $M$, there exists a Nash diffeomorphism $\phi$ from
an open semi-algebraic neighborhood $\Omega$ of the origin in
$\R^n$ onto an open semi-algebraic neighborhood $\Omega'$ of $x$
in $\R^n$ such that $\phi(0) = x$ and $\phi(\R^d \times \{0\} \cap
\Omega) = M \cap \Omega'$. }
\end{definition}
\begin{definition}
{  A \textbf{Nash map} from a Nash submanifold $M$ of $\R^m$ to a
Nash submanifold $N$ of $\R^n$ is a semi-algebraic smooth map. }
\end{definition}
\begin{remark}
Any open semi-algebraic subset of a Nash submanifold of $\R^n$ is
also a Nash submanifold of $\R^n$.
\end{remark}
\begin{theorem} \label{dimZarclos}
Let $M\subset \R^n$ be a Nash submanifold. Then it has the same
dimension as its Zarisky closure.
\end{theorem}
For proof see section 2.8 in \cite{BCR}.

Unfortunately, open semi-algebraic sets in $\R^n$ do not form a
topology, since their infinite unions are not always
semi-algebraic. This is why we need restricted topology .
\begin{definition}
{  A \textbf{$\R$-space} is a pair $(M,{\cal O}_M)$ where $M$ is a
restricted topological space and ${\cal O}_M$ a sheaf of
$\R$-algebras over $M$ which is a subsheaf of the sheaf $\R[M]$ of
real-valued functions on $M$.

A \textbf{morphism between $\R$-spaces} $(M,{\cal O}_M)$ and
$(N,{\cal O}_N)$ is a continuous map $f:M \rightarrow N$, such
that the induced morphism of sheaves $f^*:f^*(\R[N]) \rightarrow
\R[M]$ maps ${\cal O}_N$ to ${\cal O}_M$.}
\end{definition}
\begin{example}
Take for $M$ a Nash submanifold of $\R^n$, and for $ \overset
{_{\circ}} {\mathfrak S} (M)$ the family of all open subsets of
$M$ which are semi-algebraic in $\R^n$. For any open
(semi-algebraic) subset $U$ of $M$ we take as ${\cal O}_M(U)$ the
algebra ${\cal N}(U) $ of Nash functions $U \rightarrow \R$.
\end{example}
\begin{definition}
{  An \textbf {affine Nash manifold} is an $\R$-space which is
isomorphic to an $\R$-space of a closed Nash submanifold of
$\R^n$. A morphism between two affine Nash manifolds is a morphism
of $\R$-spaces between them.}
\end{definition}
\begin{example}
Any real nonsingular affine algebraic variety has a natural
structure of an affine Nash manifold.
\end{example}
\begin{remark}
Let $M \subset \R^m$ and $N \subset \R^n$ be Nash submanifolds.
Then a Nash map between them is the same as a morphism of affine
Nash manifolds between them.

Let $f:M \rightarrow N$ be a Nash map. Since an inverse of a
semi-algebraic map is semi-algebraic, $f$ is a diffeomorphism if
and only if it is an isomorphism of affine Nash manifolds.
Therefore we will call such $f$ a Nash diffeomorphism.
\end{remark}
In \cite{Shi} there is another but equivalent definition of affine
Nash manifold.
\begin{definition}
{  An \textbf {affine $C^{\infty}$ Nash manifold} is an $\R$-space
over $\R$ which is isomorphic to an $\R$-space of a Nash
submanifold of $\R^n$.}
\end{definition}
The equivalence of the definitions follows from the following
theorem.
\begin{theorem}
Any affine $C^{\infty}$ Nash manifold is Nash diffeomorphic to a
union of finite number of connected components of a real
nonsingular affine algebraic variety.
\end{theorem}
This theorem is an immediate corollary of theorem 8.4.6 in
\cite{BCR} and Theorem \ref{fincomp}.

\begin{remark}
\cite{Shi} usually uses the notion of affine $C^{\omega}$ Nash
manifold instead of affine $C^{\infty}$ Nash manifold. The two
notions are equivalent by the theorem of Malgrange (see \cite{Mal}
or Corollary I.5.7 in \cite{Shi}) and hence equivalent to what we
call just affine Nash manifold.
\end{remark}

\begin{definition}
{  A \textbf {Nash manifold} is an $\R$-space $(M,{\cal N}_M)$
which has a finite cover $(M_i)$ by open sets $M_i$ such that the
$\R$-spaces $(M_i,{{\cal N}_M}|_{M_i})$ are isomorphic to
$\R$-spaces of affine Nash manifolds.

A \textbf{morphism between Nash manifolds} is a morphism of
$\R$-spaces between them. Such morphisms are called Nash maps, and
isomorphisms are called Nash diffeomorphisms.}
\end{definition}
\begin{remark}
By Proposition \ref{fincomp}, any Nash manifold is a union of a
finite number of connected components.
\end{remark}
\begin{definition}
{  A Nash manifold is called \textbf{separated} if its restricted
topological space satisfies the standard Hausdorff separation
axiom.}
\end{definition}
\begin{remark}
Any Nash manifold has a natural structure of a smooth manifold,
and any separated Nash manifold is separated as a smooth manifold.
\end{remark}
\begin{remark}
There is a theorem by B.Malgrange (see \cite{Mal}) saying that any
Nash manifold has a natural structure of a real analytic manifold
and any Nash map between Nash manifolds is analytic. The proof is
also written on page 44 in \cite{Shi} (corollary I.5.7).
\end{remark}
\begin{example}
Any real nonsingular algebraic variety has a natural structure of
a Nash manifold.
\end{example}
\begin{proposition}
Any Nash submanifold of the projective space $\mathbb{P}^n$ is
affine.
\end{proposition}
\emph{Proof.} \\
It is enough to show that $\mathbb{P}^n$ is affine. This is
written on page 72 of \cite{BCR} (theorem 3.4.4) \proofend
\begin{remark}
So, quasiprojective Nash manifold is the same as affine Nash
manifold.
\end{remark}
%
%
\begin{notation}
{  By \textbf{open semi-algebraic subset} of a Nash manifold we
mean its open subset in the restricted topology.}
\end{notation}
The following theorem is a version of Hironaka's theorem for Nash
manifolds.
\begin{theorem}[Hironaka] \label{Hironaka}
Let $M$ be an affine Nash manifold. Then there exists a compact
affine nonsingular algebraic variety $N$ and a closed algebraic
subvariety $Z$ of $N$, which is empty if $M$ is compact, such that
$Z$ has only normal crossings in $N$ and $M$ is Nash diffeomorphic
to a union of connected components of $N - Z$.
\end{theorem}
The proof is written on page 49 of \cite{Shi} (Corollary I.5.11).
This is a consequence of Hironaka desingularization Theorem \cite{Hir}. \\
It implies the following interesting theorem.
\begin{theorem}[Local triviality of Nash manifolds]
\label{loctriv}
Any Nash manifold can be covered by finite number of open
submanifolds Nash diffeomorphic to $\R^n$.
\end{theorem}
The proof is written on page 50 of \cite{Shi} ( theorem I.5.12)
\subsection{Nash algebraic topology} \label{algtop}
In this section  we repeat known notions and theorems from
algebraic topology for the Nash case. Part of them can be found in
\cite{BCR} and \cite{Shi}.

\begin{definition}
{  Let $f:M \rightarrow N$ be a Nash map of Nash manifolds. It is
called a \textbf{Nash locally trivial fibration} with fiber $Z$ if
$Z$ is a Nash manifold and there exist a finite cover $N = \bigcup
U_i$ of $N$ by open (semi-algebraic) sets and Nash diffeomorphisms
$\phi_i$ of $f^{-1}(U_i)$ with $U_i \times Z$ such that the
composition $f \circ \phi_i^{-1}$ is the natural projection.}
\end{definition}
\begin{definition}
{  A \textbf{Nash vector bundle} over a Nash manifold $M$ is a
linear space object in the category of locally trivial fibrations
over $M$. In other words, it is an $\R$-vector bundle such that
the total space, the projection, the fiber and the trivializations
are Nash.  }
\end{definition}
\begin{remark}
In some books, for example \cite{BCR}, such vector bundles are
called pre-Nash vector bundles. They are called Nash if they can
be embedded to a trivial bundle.
\end{remark}
\begin{remark}
Direct sum and tensor product of Nash vector bundles have
canonical structure of Nash vector bundles.
\end{remark}
\begin{definition}
Let $M$ be a Nash manifold. To any locally constant sheaf
$\mathcal{F}$ on $M$ there corresponds a canonical bundle
$B(\mathcal{F})$ on $M$. Let us give an explicit construction.

Choose a cover $M = \bigcup \limits _{i=1}^k U_i$ such that
$\mathcal{F}|_{U_i}$ is isomorphic to the constant sheaf on $U_i$
with fiber $V_i$, where $V_i$ are some linear spaces. Define $N=
\bigsqcup U_i\times V_i$. We define equivalence relation: Let $u_1
\in U_i, u_2 \in U_j, v_1 \in V_i, v_2 \in V_j$ we say that
$(u_1,v_1) \sim (u_2,v_2)$ if $u_1$ and $u_2$ are the same point
in $M$ and $res_{U_i,U_i \cap U_j}(v_1) = res_{U_j,U_i \cap
U_j}(v_2)$. We define $B(F) = N/\sim$ with the obvious structure
of Nash bundle over $M$. It is easy to see that the definition
does not depend on the choice of the cover and the
trivializations.
\end{definition}
\begin{definition}
{  Let $\pi:F \rightarrow M$ be a Nash locally trivial fibration.
Consider the constant sheaf in the usual topology ${\R}^{us}_F$ on
$F$. Let $\pi_*$ denote the push functor from the category of
sheaves on $F$ to the category of sheaves on $M$.

Let $R^i \pi_*$ denote the i-th right derived functor of $\pi_*$.
Consider $R^i \pi_* (\R_M^{us})$ and restrict it to restricted
topology. We get a locally constant sheaf in the restricted
topology on $M$. We denote it by $\mathcal{H}^i(F \rightarrow
M)$.}
\end{definition}
\begin{definition}
{  Let $\pi:F \rightarrow M$ be a Nash locally trivial fibration.
Consider the constant sheaf in the usual topology $\R_F^{us}$ on
$F$. Let $\pi_!$ denote the functor of push with compact support
from the category of sheaves on $F$ to the category of sheaves on
$M$.

Consider $R^i \pi_! (\R_F^{us})$ and restrict it to restricted
topology. We get a locally constant sheaf in the restricted
topology on $M$. We denote it by $\mathcal{H}_c^i(F \rightarrow
M)$}.
\end{definition}
\begin{definition}
{  Let $\pi:F \rightarrow M$ be a Nash locally trivial fibration.
Consider the locally constant sheaf of orientations ${\mathcal
Orient}^{us}_F$ on $F$ in the usual topology. Consider the sheaf
of relative orientations ${\mathcal Orient}^{us}_F \otimes
\pi^*({\mathcal Orient}^{us}_M)$. Consider $R^i \pi_! ({\mathcal
Orient}^{us}_F \otimes \pi^*({\mathcal Orient}^{us}_M))$ and
restrict it to restricted topology. We get a locally constant
sheaf in the restricted topology on $M$. We denote it by
$\mathcal{TH}_c^i(F \rightarrow M)$.}
\end{definition}
\begin{notation}
{  We denote $H^i(F \rightarrow M) := B(\mathcal{H}^i(F
\rightarrow M))$, $H^i_c(F \rightarrow M) := B(\mathcal{H}^i_c(F
\rightarrow M))$, $TH_c^i(F \rightarrow M) := B(\mathcal{TH}_c^i(F
\rightarrow M))).$}
\end{notation}
\begin{proposition}
Tangent, normal and conormal bundles, the bundle of differential
$k$-forms, the bundle of orientations, the bundle of densities,
etc. have canonical structure of Nash bundles.
\end{proposition}
For proof see e.g. Theorems 3.4.3 and 3.4.4 in \cite{AG}.
\begin{notation}
{  Let $M$ be a Nash manifold. We denote by \textbf{$ Orient_M$}
the bundle of orientations on $M$ and by \textbf{$ D_M$} the
bundle of densities on $M$.}
\end{notation}

\begin{notation}
{  Let $E\rightarrow M$ be a Nash bundle. We denote
\textbf{$\widetilde{E}$}$:=E^*\otimes D_M$.}
\end{notation}

%
Using Hironaka theorem (\ref{Hironaka}) we will prove the
following useful result.
\begin{theorem}
Let $M$ be a Nash manifold. Then
$H^i(M)$,$H_c^i(M)$,$H_c^i(M,{\mathcal Orient}^{us}_M)$ are finite
dimensional.
\end{theorem}
For this we will need the following lemma.
\begin{lemma}
Let $M$ be a smooth manifold. Let $N \subset M$ be a closed submanifold and denote 
$U=M - N$. Let $\mathcal{F}$ be a locally constant sheaf on $M$
such that $H_c^i(M,\mathcal{F})$ and $H_c^i(N,\mathcal{F}|_N)$ are
finite dimensional, where by $\mathcal{F}|_N$ we mean restriction
as a local system. Then $H_c^i(U,\mathcal{F}|_U)$ is also finite
dimensional.
\end{lemma}
\emph{Proof of the lemma}. Let $\phi:N \rightarrow M$ and $\psi:U
\rightarrow M$ be the standard imbeddings. Note that $\phi_!$ and
$\psi_!$ are exact and $\phi_*=\phi_!$ so $H_c^i(U,\mathcal{F}|_U)
\cong H_c^i(M,j_!(\mathcal{F}|_U))$ and
$H_c^i(N,\mathcal{F}|_N)\cong H_c^i(M,i_!(\mathcal{F}|_N))$. So
from short exact sequence $0\rightarrow
\psi_!(\mathcal{F}|_U)\rightarrow \mathcal{F}\rightarrow
\phi_!(\mathcal{F}|_N)\rightarrow 0$ of sheaves on $M$ we see that
$H_c^i(U,\mathcal{F}|_U)$ is also finite dimensional. \proofend\\
\emph{Proof of the theorem}. Intersection of affine open Nash
submnaifolds is affine, hence by Mayer - Vietories long exact
sequence (see e.g. \cite{BT}, section I.2) it is enough to prove
the theorem for affine Nash manifolds. Note that
$H_c^i(M)=H_c^i(M,{\R}_M^{us})$ where ${\R}_M^{us}$ is the
constant sheaf on $M$. Now using the lemma and Hironaka theorem we
can easily show by induction that $H_c^i(M)$ and
$H_c^i(M,{\mathcal Orient}^{us}_M)$ are finite dimensional. By
Poincar\'{e} duality $H^i(M) \cong H_c^i(M,{\mathcal
Orient}^{us}_M)^*$ and hence is
finite dimensional. \proofend \\
Now we will give another definition of Nash
locally trivial fibration.
\begin{definition}
{  Let $f:M \rightarrow N$ be a Nash map of Nash manifolds. It is
called a \textbf{Nash locally trivial fibration} if there exist a
Nash manifold $M$ and surjective submersive Nash map $g:M
\rightarrow N$ such that the base change $h:M \underset
{N}{\times} M \rightarrow M$ is trivializable, i.e. there exists a
Nash manifold $Z$ and a Nash diffeomorphism $k:M \underset
{N}{\times} M \rightarrow M \times Z$ such that $\pi \circ k = h$
where $\pi : M \times Z \rightarrow M$ is the standard
projection.}
\end{definition}
In order to prove that this definition is equivalent to the
previous one, it is enough to prove the following theorem.
\begin{theorem} \label{sursec}
Let $M$ and $N$ be Nash manifolds and $\nu:M \rightarrow N$ be a
surjective submersive Nash map. Then locally it has a Nash
section, i.e. there exists a finite open cover $N= \bigcup \limits
_{i=1}^k U_i$ such that $\nu$ has a Nash section on each $U_i$.
\end{theorem}
For proof see Appendix \ref{AppSurSec}.

\subsection{Schwartz functions on Nash manifolds} \label{Schwartz}

In this section we will review some preliminaries on Schwartz
functions on Nash manifolds defined in \cite{AG}.

The \Fre space ${\Sc}(\R^n)$ of Schwartz functions on $\R^n$ was
defined by Laurant Schwartz to be the space of all smooth
functions that decay faster than $1/|x|^n$ for all n.

In \cite{AG} we have defined for any Nash manifold $M$ the \Fre
space ${\Sc}(M)$ of Schwartz functions on it.

As Schwartz functions cannot be restricted to open subsets, but
can be continued by 0 from open subsets, they form a cosheaf
rather than a sheaf.

We have defined for any Nash bundle $E$ over $M$ the cosheaf
${\Sc}_M^E$ over $M$ (in the restricted topology) of Schwartz
sections of $E$. These cosheaves satisfy:
${\Sc}_M^E(U)={\Sc}(U,E|_U) = \{\xi \in {\Sc}(M,E) | \xi$ vanishes
with all its derivatives on $M - U \}$. We have also defined the
sheaf $\G_M^E$ of generalized Schwartz sections of $E$ by
$\G_M^E(U)= ({\Sc}_M^{\widetilde{E}}(U))^*$. This sheaf is flabby.

The fact that ${\Sc}_M^E$  satisfies the cosheaf axioms follows
from the following version of partition of unity:

\begin{theorem}[Partition of unity] \label {partuni}
Let $M$ be a Nash manifold, and $(U_i)_{i=1}^n$ - finite open
 cover by affine Nash submanifolds. Then there exist smooth functions
$\alpha_1,...,\alpha_n$ such that $supp(\alpha_i) \subset U_i$,
$\sum \limits _{i=1}^n \alpha_i = 1$ and for any $g \in {\Sc}(M)$,
$\alpha_ig \in {\Sc}(U_i)$.
\end{theorem}
For proof see \cite{AG}, section 5.2. \\
We will use the following proposition which follows trivially from
the definition of the sheaves of Schwartz sections and generalized
Schwartz sections given in \cite{AG}, section 5.
\begin{proposition}
Let $\mathcal{F}$ be a locally constant sheaf over a Nash manifold
$M$. Then ${\Sc}_M^{B(\mathcal{F})} \cong {\Sc}_M \otimes
\mathcal{F}'$ and $\G_M^{B(\mathcal{F})} \cong \G_M \otimes
\mathcal{F}$. Moreover, if $E$ is a Nash vector bundle over $M$
then ${\Sc}_M^{B(\mathcal{F})\otimes E} \cong {\Sc}_M^E \otimes
\mathcal{F}'$ and $\G_M^{B(\mathcal{F})\otimes E} \cong \G_M^E
\otimes \mathcal{F}$. Recall that $B(\mathcal{F})$ is the bundle
corresponding to $\mathcal{F}$ and $\mathcal{F'}$ is the cosheaf
corresponding to $\mathcal{F}$ as they were defined in sections
\ref{algtop} and \ref{sheaftheo}.
\end{proposition}

To conclude, we will list the important statements from \cite
{AG}:

\begin{property} \label{p1}
Compatibility: For open semi-algebraic subset $U \subset M$,
${\Sc}^E_M(U)={\Sc}(U,E|_U)$.
\end{property}

\begin{property}\label{p2}
$\Sc(\R ^n)$ = Classical Schwartz functions on $\R ^n$.
\end{property}

\begin{property} \label{p3}
For compact $M$, ${\Sc}(M,E) = $ smooth global sections of $E$.
\end{property}

\begin{property} \label{p4}
$\G_M^{E}=({\Sc}^{\widetilde {E} }_M)^*$ , where $\widetilde
{E}=E^*\otimes D_M$ and $D_M$ is the bundle of densities on M .
\end{property}

\begin{property} \label{p5}
Let $Z \subset M$ be a Nash closed submanifold. Then restriction
maps  ${\Sc}(M,E)$ onto ${\Sc}(Z,E|_Z)$.
\end{property}

\begin{property} \label{p6}
Let $U \subset M$  be a semi-algebraic open subset. Then
$${\Sc}^E_M(U) \cong \{\xi \in {\Sc}(M,E)| \quad \xi \text{ is 0 on } M - U \text{ with all derivatives} \}.$$
\end{property}

\begin{property} \label{p7}
Let $Z \subset M$ be a Nash closed submanifold. Consider $V=\{\xi
\in \G(M,E) |\xi$ is supported in $Z \}$. It has canonical
filtration $V_i$ such that its factors are canonically isomorphic
to $\G(Z,E|_Z \otimes {Sym}^i(CN_Z^M) \otimes {D_M^*|}_Z \otimes
D_Z)$ where $CN_Z^M$ is the conormal bundle of $Z$ in $M$ and
$Sym^i$ means i-th symmetric power.
\end{property}
\subsection{Nuclear \Fre spaces} \label{NucFre}

\begin{definition}
We call a complex of topological vector spaces \textbf{admissible}
if all its differentials have closed images.
\end{definition}

We will need the following classical facts from the theory of
nuclear \Fre spaces.

\begin{itemize}
\item Let $V$  be a nuclear \Fre space and $W$ be a closed subspace. Then both $W$ and $V/W$ are nuclear \Fre
spaces.

\item Let
$$\mathcal{C}: 0 \rightarrow
C_1 \rightarrow ... \rightarrow C_n \rightarrow 0$$
be an admissible complex of nuclear \Fre spaces. Then the complex
$\mathcal{C}^*$ is also admissible and $H^i(\mathcal{C}^* ) \cong
H^i( \mathcal{C}_i) ^*.$

\item Let $V$  be a nuclear \Fre space. Then the complex $ \mathcal{C} \ctp V $ is an admissible
complex of nuclear \Fre spaces and $H^i( \mathcal{C} \otimes V)
\cong H^i( \mathcal{C}) \otimes V .$

\item $\Sc(\R^n)$ is a nuclear \Fre space.

\item $\Sc(\R^{n+m})= \Sc(\R^n) \ctp \Sc(\R^m)$.
\end{itemize}

A good exposition on nuclear \Fre spaces can be found in Appendix
A in \cite{CHM}.

\begin{corollary}
Let $M$ be a Nash manifold and $E$ be a Nash bundle over it. Then
$\Sc(M,E)$ is a nuclear \Fre space.
\end{corollary}
\emph{Proof.} By definition of $\Sc(M,E)$ and by Theorem
\ref{loctriv}, $\Sc(M,E)$ is a quotient of direct sum of several
copies of $\Sc(\R^n)$. \proofend

\begin{corollary}
Let $M_i$, $i=1,2$ be Nash manifolds and $E_i$ be Nash bundles
over $M_i$. Then $$\Sc(M_1 \times M_2,E_1 \boxtimes
E_2)=\Sc(M_1,E_1) \ctp \Sc(M_2,E_2),$$ where $E_1 \boxtimes E_2$
denotes the exterior product.
\end{corollary}

\section{De-Rham theorem for Schwartz functions on Nash
manifolds} \label{DeRham}

\subsection{De-Rham theorem for Schwartz functions on Nash
manifolds}
\begin{theorem} \label{quasiiso}
Let $M$ be an affine Nash manifold. Consider the de-Rham complex
of $M$ with compactly supported coefficients
$$DR_c(M): 0 \rightarrow
C^{\infty}_c(M,{\Omega^0_M}) \rightarrow ... \rightarrow
C^{\infty}_c(M,{\Omega^n_M}) \rightarrow 0$$ and the natural map
$i: DR_c(M) \rightarrow DR_{\Sc}(M)$.  Then $i$ is a
quasiisomorphism, i.e. it induces an isomorphism on the
cohomologies.
\end{theorem}
\emph{Proof.} Let $N\supset M$ be the compactification of $M$
given by Hironaka theorem, i.e. $N$ is a compact Nash manifold,
$N$ = $M \udot D \udot U$ where $M$ and $U$ are open and $D =
\bigcup \limits _{i=1}^k D_i$ where $D_i\subset N$ is a closed
Nash submanifold of codimension 1 and all the intersections are
normal, i.e. every $y\in N$  has a neighborhood $V$ with a
diffeomorphism $\phi:V \rightarrow \R^n$ such that $\phi(D_i\cap
V)$ is either a coordinate hyperplane or empty. Denote $Z=N- M$.

$N$ has a structure of compact smooth manifold. We build two
complexes $\mathcal{DR}_1$ and $\mathcal{DR}_2$ of sheaves on $N$
in the classical topology by $\mathcal{DR}_1^k(W):=\{\omega \in
C^{\infty}(W,\Omega^k)|\omega$ vanishes in a neighborhood of $Z\}$
and $\mathcal{DR}_2^k(W)=\{\omega \in
C^{\infty}(W,\Omega^k)|\omega$ vanishes on $Z$ with all its
derivatives $\}$. As the differential we take the standard de-Rham
differential.

Note that we have a natural embedding of complexes ${\mathcal
I}:\mathcal{DR}_1 \rightarrow \mathcal{DR}_2$. Note also that
$\mathcal{DR}_1(M) \cong DR_c(M)$ and $\mathcal{DR}_2(M) \cong
DR_{\Sc}(M)$. The theorem follows from the facts that
$\mathcal{DR}_{1,2}^i$ are $\Gamma$ - acyclic sheaves and that
${\mathcal I}$ is a quasiisomorphism. Let us prove these two facts
now.

$\mathcal{DR}_{1,2}^i$ are fine (i.e. have partition of unity),
which follows from the classical partition of unity. So, by
theorem 5.25 from \cite{War} they are acyclic.

The statement that ${\mathcal I}$ is a quasiisomorphism is a local
statement, so we will verify that ${\mathcal I}:
\mathcal{DR}_{1}(W) \to \mathcal{DR}_{2}(W)$ is a quasiisomorphism
for small enough $W$. Since all the intersections in $D$ are
normal, it is enough to check it for the case $W \cong \R^n$ and
$D \cap W$ is a union of coordinate hyperplanes. In this case, the
proof is technical and all its ideas are taken from classical
proof of Poincar\'{e} lemma. We will give it now only for
completeness and we recommend the reader to
skip to the end of the proof. \\
$(N- D) \cap W$ splits to a union of connected components of the
form $R_{>0}^k \times \R^l$. Hence complexes
$\mathcal{DR}_{1,2}(W)$ split to direct sum of the complexes
corresponding to the connected
components. Therefore, it is enough to check this statement in the following two cases: \\
Case 1 $W=R_{\geqslant 0}^k \times \R^l$, $U \cap W= R_{>0}^k
\times
\R^l$, $M \cap W = \emptyset$ \\
Case 2 $W=R_{\geqslant 0}^k \times \R^l$, $U \cap W=
\emptyset$, $M \cap W= R_{>0}^k \times \R^l$.\\
Case 1 is trivial, as $\mathcal{DR}_1(W) = \mathcal{DR}_2(W) = 0$ in this case. \\
Case 2: If $k=0$ then ${\mathcal I}=Id$. Otherwise we will show
that the cohomologies of both complexes vanish. Clearly $H^0_{1,2}
= 0$ since the only constant function which vanishes on $D$ is 0.
Now let $\omega \in \mathcal{DR}_{1,2}^m(W)$ be a closed form. We
can write $\omega$ in coordinates $dx_1,...,dx_{k+l}$: $\omega =
\omega_1 \wedge dx_1 + \omega_2$ where neither $\omega_1$ nor
$\omega_2$ contain $dx_1$. Let $f_j$ be the coefficients of
$\omega_1$ and define $g_j(x_1,...,x_{k+l})= \int \limits _0
^{x_1} f_j(t,x_2,...,x_{k+l})dt$ and let $\lambda$ be the form
with coefficients $g_j$. It is easy to check that $d\lambda =
\omega$ and $\lambda \in \mathcal{DR}_{1,2}^{m-1}(W)$. \proofend

\begin{theorem} \label{quasiisogen}
Let $M$ be an affine Nash manifold. Consider the de-Rham complex
of $M$ with coefficients in classical generalized functions, i.e.
functionals on compactly supported densities.
$$DR_{-\infty}(M): 0 \rightarrow
C^{-\infty}(M,{\Omega^0_M}) \rightarrow ... \rightarrow
C^{-\infty}(M,{\Omega^n_M}) \rightarrow 0$$ and the natural map
$i: DR_{\G}(M) \to DR_{-\infty}(M)$.  Then $i$ is a
quasiisomorphism.
\end{theorem}

\noindent \emph{Proof.} Let $N$, $D$, $D_i$, $U$ and $Z$ be the
same as in the proof of Theorem \ref{quasiiso}. We again build two
complexes $\mathcal{DR}_1$ and $\mathcal{DR}_2$ of sheaves on $N$
in the classical topology by $\mathcal{DR}_1^k(W):=k$-forms on $W-
Z$ with generalized coefficients and
$\mathcal{DR}_2^k(W):=k$-forms with generalized coefficients on
$W$ modulo $k$-forms with generalized coefficients on $W$
supported in $Z\cap W$. We have an embedding ${\mathcal
I}:\mathcal{DR}_2 \hookrightarrow \mathcal{DR}_1$. Again, by
classical partition of unity the sheaves are fine and hence
acyclic, so it is enough to prove that ${\mathcal I}$ is a
quasiisomorphism. Again, we check it locally and the only
interesting case is $W=R_{\geq 0}^k \times \R^l$, $U \cap W=
\emptyset$, $M \cap W = R_{>0}^k \times \R^l$ where $k>0$. Define
a map $\phi :\R \rightarrow \mathcal{DR}(W)^0_{1,2}$ by setting
$\phi(c)$ to be the constant generalized function $c$.  It gives
us extensions $\widetilde{\mathcal{DR}_{1,2}}(W)$ of complexes
$\mathcal{DR}_{1,2}(W)$ and $\widetilde{{\mathcal I}}$ of
${\mathcal I}$. It is enough to prove that $\widetilde{{\mathcal
I}}$ is a quasiisomorphism. For this we will prove that both
complexes are acyclic. Fix standard orientation on $N$. Now our
complexes become dual to $$C_1: 0 \leftarrow \R \leftarrow
C^{\infty}_c(W\cap M,{\Omega^n_{W\cap M}}) \leftarrow ...
\leftarrow C^{\infty}_c(W\cap M,{\Omega^0_{W\cap M}}) \leftarrow
0$$ and
$$C_2: 0 \leftarrow \R \leftarrow C^{\infty}_c(W,W\cap D,{\Omega^n_W})
\leftarrow ... \leftarrow C^{\infty}_c(W,W\cap
D,{\Omega^0_W})\leftarrow 0 $$ where $C^{\infty}_c(W,W\cap
D,{\Omega^n_W})$ are compactly supported forms which vanish with
all their derivatives on $W\cap D$. We will prove that $C_{1,2}$
are homotopically equivalent to zero and this will give us that
$\widetilde{\mathcal{DR}_{1,2}}(W)$ are also homotopically
equivalent to zero and hence are acyclic. The complex $C_1$ is
isomorphic to the following complex
$$ C'_1: 0 \leftarrow \R \leftarrow
C^{\infty}_c({\R^n},{\Omega^n_{\R^n}}) \leftarrow ... \leftarrow
C^{\infty}_c({\R^n},{\Omega^0_{\R^n}}) \leftarrow 0 .$$ In section
I.4 of \cite{BT} (Poincar\'{e} lemma for compactly supported
cohomologies) it is proven that $C_1$ is homotopy equivalent to
zero. In the same way we
can prove that $C_2$ is homotopically equivalent to zero.

\proofend

The following theorem is classical.
\begin{theorem} \label{quasiisoclas}
Let $M$ be a smooth manifold. Consider the de-Rham complex of $M$
with coefficients in classical generalized functions
$DR_{-\infty}(M)$, the de-Rham complex of $M$ with coefficients in
smooth functions $DR(M)$ and the natural map $i: DR(M) \rightarrow
DR_{-\infty}(M)$. Then $i$ is a quasiisomorphism.
\end{theorem}
\emph{Proof. }Let $\mathcal{DR}_{-\infty}$ and $\mathcal{DR}$ be
the de-Rham complex of $M$ with coefficients in the sheaves of
classical generalized functions and smooth functions
correspondingly. The sheaves in these complexes are acyclic hence
it is enough to show that the natural map ${\mathcal I}:
\mathcal{DR} \rightarrow \mathcal{DR}_{-\infty}$ is a
quasiisomorphism. This is proven by a local computation similar to
the one in the proof of the last theorem. \proofend
\begin{definition} 
Let $M$ be a Nash manifold. We define the \textbf{twisted bundle
of $k$-differential forms on $M$ by $T
\Omega^k_M:=\Omega^k_M\otimes Orient_M$} and correspondingly the
\textbf{twisted de-Rham complexes} $$TDR_{-\infty}(M), \,
TDR_{\G}(M), \, TDR(M), \, TDR_{{\Sc}}(M), \, TDR_{c}(M).$$
\end{definition}
\begin{remark}
{  Note that $T \Omega^{n-k}_M \cong \widetilde{\Omega^{k}_M}$.
This gives us a natural pairing between ${\Sc}(M,T
\Omega^{n-k}_M)$ and $\G(M,\Omega^k_M)$}.
\end{remark}
\begin{remark}
The theorems \ref{quasiiso}, \ref{quasiisogen} and
\ref{quasiisoclas} hold true also for the twisted de-Rham
complexes and the proofs are the same.
\end{remark}
The bottom line of this section is the following version of
de-Rham theorem
\begin{theorem} \label{bottomline}
Let $M$ be an affine Nash manifold of dimension $n$. Then
\\ $H^i(DR_{\G}(M))\cong H^i(M)$
\\ $H^i(DR_{{\Sc}}(M))\cong H_c^i(M)$
\\ $H^i(TDR_{{\Sc}}(M))\cong  H_c^i(M,{\mathcal Orient}^{us}_M)$
\\and the pairing between
$\G(M,{\Omega_M^i})$ and ${\Sc}(M,T \Omega_M^{n-i})$ gives an
isomorphism between $H^i(DR_{\G}(M))$ and
$(H^{n-i}(TDR_{{\Sc}}(M)))^*$.
\end{theorem}

\noindent \emph{Proof.} The theorem is a direct corollary from
theorems \ref{quasiiso} \ref{quasiisogen} \ref{quasiisoclas}  for
the standard and the twisted cases and from classical Poincar\'{e}
duality. \proofend

\begin{corollary}
The complexes $DR_{\G}(M), DR_{{\Sc}}(M)$ and $TDR_{{\Sc}}(M)$ are
admissible.
\end{corollary}

\subsection{Relative de-Rham theorem for Nash locally trivial fibration}
\begin{definition} 
Let $F \overset{\pi}{\rightarrow} M$ be a locally trivial
fibration. Let $ E \rightarrow M$ be a Nash bundle. We define
$T_{F \rightarrow M} \subset T_F$ by $T_{F \rightarrow M}=ker(d
\pi)$. We denote %
$$\Omega^{i,E}_{F \rightarrow M}:=((T_{F \rightarrow M})^*)^{\wedge i}\otimes \pi^*E , \,Orient_{F
\rightarrow M}=Orient_F\otimes \pi^*(Orient_M) \text{ and }T
\Omega^{i,E}_{F \rightarrow M}:=\Omega^{i,E}_{F \rightarrow
M}\otimes Orient_{F \rightarrow M}.$$
Now we can define \textbf{the relative de-Rham complexes}
$$DR^E_{\G}(F \rightarrow M),\, DR^E_{{\Sc}}(F \rightarrow M), \,TDR^E_{\G}(F \rightarrow
M), \, TDR^E_{{\Sc}}(F \rightarrow M).$$
If $E$ is trivial we will omit it.
\end{definition}

The goal of this section is to prove the following theorem.
\begin{theorem} \label{RelDeRham}
Let $p:F \rightarrow M$ be a Nash locally trivial fibration. Then
\setcounter{equation}{0}
\begin{align}\label{i}
&H^k(DR_{\Sc}^E(F\rightarrow M)) \cong {\Sc}(M,H^k_c(F\rightarrow
M)\otimes E).\\ \label{ii} & H^k(TDR_{\Sc}^E(F\rightarrow M))
\cong {\Sc}(M,TH_c^k(F\rightarrow M)\otimes E).\\
 \label{iii} &H^k(DR_{\G}^E(F\rightarrow M)) \cong \G(M,H^k(F\rightarrow
M)\otimes E).
\end{align}
\end{theorem}

\noindent \emph{Proof}

\noindent (\ref{i}) Step 1. Proof for the case $M=\R^n$, the
fibration $F \to M$ is trivial and $E$ is trivial.

It follows from Theorem \ref{bottomline} using subsection
\ref{NucFre}.

Step 2. Proof in the general case.

Let $C_i \subset \Sc(F, \Omega_{F \to M}^{i,E})$ be the subspace
of closed forms. We have to construct a continuous onto map
$\phi_i:C_i \twoheadrightarrow \Sc(M, H^i(F \to M) \otimes E)$
whose kernel is the space of exact forms. Fix a cover $M =
\bigcup_{k=1}^m U_k$ such that $U_k$ are Nash diffeomorphic to
$\R^n$ and $F|_{U_k}$ and $E|_{U_k}$ are trivial. Fix a partition
of unity $1 = \sum \alpha_i$ such that and for any $g \in
{\Sc}(F)$, $\alpha_ig \in {\Sc}(p^{-1}(U_i))$. Note that for any
$\omega \in \Sc(F,\Omega_{F \to M}^{i,E})$ we have $\alpha_i
\omega \in \Sc(p^{-1}(U_i),\Omega_{F \to M}^{i,E})$.  By the
previous step,
$$ H^i(DR^{E|_{U_k}}_{F|_{U_k} \to U_k}) \cong \Sc(U_k,H^i(F|_{U_k} \to U_k)
\otimes E|_{U_k}).$$ 
For any form $\nu \in \Sc(p^{-1}(U_k),\Omega_{F|_{U_k} \to
U_k}^{i,E|_{U_k}})$ we consider the class $[\nu]$ as an element
$$[\nu] \in \Sc(U_k,H^i(F|_{U_k} \to U_k) \otimes E|_{U_k}) \subset \Sc(M,
H^i(F \to M) \otimes E).$$

Now let $\omega \in C_i$. Define $$ \phi_i(\omega):=\sum_{k=1}^m
[\alpha_i \omega].$$ It is easy to see that $\phi$ satisfies the
requirements and does not depend on the choice of $U_k$ and
$\alpha_k$.

\noindent (\ref{ii}) Is proven in the same as (\ref{i}).

\noindent (\ref{iii}) follows from (\ref{ii}) using subsection
\ref{NucFre}.
\section{Shapiro lemma} \label{secShapLem}
\setcounter{lemma}{0} In this section we formulate and prove a
version of Shapiro lemma for generalized Schwartz sections of Nash
equivariant bundles.

\begin{definition} 
{  Let $\mathfrak{g}$ be a Lie algebra of dimension $n$. Let
$\rho$ be its representation. Define $H^i(\mathfrak{g},\rho)$ to
be the cohomologies of the complex: $$C(\mathfrak{g},\rho): 0
 {\rightarrow} \rho  {\rightarrow}
\mathfrak{g}^* \otimes \rho  {\rightarrow} ({\mathfrak{g}^*})
^{\wedge 2} \otimes \rho  {\rightarrow}...  {\rightarrow}
({\mathfrak{g}^*}) ^{\wedge n} \otimes \rho  {\rightarrow} 0$$
with the differential defined by
\begin{multline*}
d\omega(x_1,...,x_{n+1}) = \sum _{i=1}^{n+1} (-1)^i \rho(x_i)
\omega(x_1,...,x_{i-1},x_{i+1},...,x_{n+1}) +\\ +\sum
_{i<j}(-1)^{i+j}\omega([x_i,x_j],x_1,...,x_{i-1},x_{i+1},...,x_{j-1},x_{j+1},...,x_{n+1})
\end{multline*} where we interpret $({\mathfrak{g}}^*) ^{\wedge k}
\otimes \rho$ as anti-symmetric $\rho$-valued $k$-forms on
${\mathfrak{g}}$.}
\end{definition}
\begin{remark}
$H^i(\mathfrak{g},\rho)$ is the i-th derived functor of the functor
$\rho \mapsto {\rho}^{\mathfrak{g}}$.
\end{remark}
\begin{definition}
{  A \textbf{Nash group} is a group object in the category of Nash
manifolds, i.e. a Nash manifold $G$ together with a point $e \in
G$ and Nash maps $m:G \times G \rightarrow G$ and $inv:G
\rightarrow G$ which satisfy the standard group axioms.

A \textbf{Nash $G$- manifold} is a Nash manifold $M$ together with
a Nash map $a:G \times M \rightarrow M$ satisfying $a(gh,x)=
a(g,a(h,x))$.

A \textbf{Nash $G$ - equivariant bundle} is a Nash vector bundle
$E$ over a Nash $G$-manifold $M$ together with an isomorphism of
Nash bundles $pr^*(E) \simeq a^*(E)$ where $pr:G \times M
\rightarrow M$ is the standard projection.}
\end{definition}
\begin{definition}
{  Let $G$ be a Nash group and $M$ be a Nash $G$ manifold. We
define \textbf{the quotient space $G \setminus M$ }to be the
following $\R$-space. As a set, it is the set theoretical
quotient. A subset $U \subset G \setminus M$ is open iff
$\pi^{-1}(U)$ is open, where $\pi$ is the standard projection $M
\rightarrow G \setminus M$. The sheaf of regular functions is
defined by $\mathcal{O}(U) = \{f | f \circ \pi \in
\mathcal{N}(\pi^{-1}(U))\}$.}
\end{definition}
\begin{definition}
{A Nash action of a Nash group $G$ on a Nash manifold $M$ is
called \textbf{strictly simple} if it is simple (i.e. all
stabilizers are trivial) and $G \setminus M$ is a separated Nash
manifold.}
\end{definition}
\begin{proposition}
Let $G$ be a Nash group and $M$ be a Nash $G$ manifold. Suppose
that the action is strictly simple. Then the projection $\pi:M
\rightarrow G \setminus M$ is a Nash locally trivial fibration.
\end{proposition}
\emph{Proof.} From differential topology we know that $\pi$ is a
submersion. Consider the base change $M \underset{G \setminus
M}{\times} M \rightarrow M$. It is Nash diffeomorphic to the
trivial projection $M \times G \rightarrow M$. \proofend
\begin{corollary}
Let $G$ be a Nash group and $M$ be a Nash $G$ manifold with
strictly simple action. Let $N$ be any $G$ manifold. Then the
diagonal action on $M \times N$ is strictly simple.
\end{corollary}
\emph{Proof.} If the fibration $M\rightarrow G \setminus M$ is
trivial, the statement is clear. It is locally trivial by the
proposition, and the statement is local on $G \setminus M$.
\proofend
\begin{remark}
Let $G$ be a Nash group, $M$ be a Nash $G$-manifold and $E
\rightarrow M$ be a Nash $G$-equivariant bundle. Then the spaces
$\Sc(M,E)$ and $\G(M,E)$ have natural structure of
$G$-representations. Moreover, they are smooth $G$-representations
and hence they have a natural structure of
$\mathfrak{g}$-representations where $\mathfrak{g}$ is the Lie
algebra of $G$.
\end{remark}
Now we give a recipe how to compute cohomologies of such
representations.
\begin{theorem}\label{recipe}
Let $G$ be a Nash group. Let $M$ be a Nash $G$-manifold and $E
\rightarrow M$ be a Nash $G$-equivariant bundle. Let $N$ be a
strictly simple Nash $G$-manifold. Suppose that $N$ and $G$ are
cohomologically trivial (i.e. all their cohomologies except $H^0$
vanish and $H^0=\R$) and affine. Denote $F=M\times N$ . Note that
the bundle $E \boxtimes \Omega_N^i$ has Nash $G$-equivariant
structure given by diagonal action. Hence the relative de-Rham
complex $DR_{\G}^E(F \rightarrow M)$ is a complex of
representations of $\mathfrak{g}$. Then $H^i(\mathfrak{g},
\G(M,E)) = H^i((DR_{\G}^E(F \rightarrow M))^{\mathfrak{g}})$.
\end{theorem}
For this theorem we will need the following lemma.
\begin{lemma} \label{DRgroup}
Let $G$ be a Nash group. Let $F$ be a strictly simple Nash
$G$-manifold. Denote $M:=G \setminus F$ let $E \rightarrow M$ be a
Nash bundle. Then the relative de-Rham complex
$DR_{\G}^E(F\rightarrow M)$ is isomorphic to the complex
$C(\mathfrak{g},\G(F,\pi^*E))$, where $\pi:F\rightarrow M$ is the
standard projection.
\end{lemma}
\emph{Proof.} By partition of unity it is enough to prove for the
case that the fibration $\pi:F\rightarrow M$ is trivial. In this
case we can imbed $\mathfrak{g}$ into the space of Nash sections
of the bundle $T_{F \rightarrow M}\rightarrow F$ and its image
will generate the space of all Nash sections of $T_{F \rightarrow
M}\rightarrow F$ over $\mathcal{N}(F)$. This gives us an
isomorphism between $\G(F)\otimes \mathfrak{g}$ and $\G(F,T_{F
\rightarrow M})$ and in the same way between $({\mathfrak{g}^*})
^{\wedge k} \otimes \G(\pi^*E,F)$ and $\G(F,\Omega^{k,E}_{F
\rightarrow M})$. It is easy to check that the last isomorphisms
form an isomorphism of complexes between $DR_{\G}^E(F\rightarrow
M)$ and $C(\mathfrak{g},\G(F,\pi^*E))$.

$ $ \proofend\\
\emph{Proof of Theorem \ref{recipe}}. From relative de-Rham
theorem (\ref{RelDeRham}), we know that the complex $DR_{\G}^E(F
\rightarrow M)$ is a resolution of $\G(M,E)$ (i.e. all its higher
cohomologies vanish and the 0's cohomology is equal to $\G(M,E)$).
So it is enough to prove that the representations $\G(F,E
\boxtimes \Omega_N^i)$ are $\mathfrak{g}$- acyclic. Denote $Z:=G
\setminus F$. The fact that the bundle $E \boxtimes
\Omega_N^i\rightarrow F$ is $G$- equivariant gives us an action of
$G$ on the total space $E \boxtimes \Omega_N^i$. Denote $B:= G
\setminus (E \boxtimes \Omega_N^i)$. Note that $B \rightarrow Z$
is a Nash bundle and $F \rightarrow Z$ is a Nash locally trivial
fibration. By the lemma, the complex $C(\mathfrak{g},\G(E
\boxtimes \Omega_N^i,F))$ is isomorphic to the relative de-Rham
complex $DR_{\G}^B(F \rightarrow Z)$ and again by relative de-Rham
theorem $H^i(DR_{\G}^B(F \rightarrow Z))=0$ for $i>0$. \proofend
\begin{proposition}
Let $G$ be a connected Nash group and $F$ be a Nash $G$ manifold
with strictly simple action. Denote $M:=G \setminus F$ and let $E
\rightarrow M$ be a Nash bundle. Then
$(\G(F,\pi^*(E)))^{\mathfrak{g}}\cong \G(M,E)$ where
$\pi:F\rightarrow M$ is the standard projection.
\end{proposition}
\emph{Proof.} It is a direct corollary of Lemma \ref{DRgroup} and
relative de-Rham theorem (\ref{RelDeRham})\proofend
\begin{corollary}\label{H2G}
Let $G$ be a Nash group and $M$ be a transitive Nash $G$ manifold.
Let $x \in M$ and denote $H:=stab_G(x)$. Consider the diagonal
action of $G$ on $M\times G$. Let $E\rightarrow M\times G$ be a
$G$ equivariant Nash bundle. Then $\G(M\times
G,E)^\mathfrak{g}\cong \G(\{x\} \times G,E|_{\{x\} \times
G})^\mathfrak{h}$.
\end{corollary}
Now we can prove Shapiro lemma.
\begin{theorem}[Shapiro lemma]
Let $G$ be a Nash group and $M$ be a transitive Nash $G$ manifold.
Let $x \in M$ and denote $H:=stab_G(x)$. Let $E\rightarrow M$ be a
$G$ equivariant Nash bundle. Let $V$ be the fiber of $E$ in $x$.
Suppose $G$ and $H$ are cohomologically trivial. Then
$H^i(\mathfrak{g},\G(M,E))\cong H^i(\mathfrak{h},V)$.
\end{theorem}
\emph{Proof.} From the recipe of computing cohomologies (Theorem
\ref{recipe})  we see that
$$H^i(\mathfrak{g},\G(M,E))\cong H^i((DR^E_{\G}(M\times
G\rightarrow M))^{\mathfrak{g}}).$$ By Corollary \ref{H2G}
$$H^i((DR^E_{\G}(M\times G\rightarrow M))^\mathfrak{g}) \cong
H^i((DR^V_{\G}(\{x\}\times G\rightarrow \{x\}))^{\mathfrak{h}})$$
and again by the recipe of computing cohomologies (Theorem
\ref{recipe})
$$H^i((DR^V_{\G}(\{x\}\times G\rightarrow
\{x\}))^\mathfrak{h})\cong H^i(\mathfrak{h},V).$$\proofend \\
To make the theorem complete we need to prove that a quotient of a
Nash group by its Nash subgroup is a Nash manifold. We prove it in
the case of linear Nash group.
\begin{proposition} \label{takoyto}
Let $H <G<GL_n$ be Nash groups. Then the action of $H$ on $G$ is
strictly simple.
\end{proposition}
To prove the proposition we will need the following lemma.
\begin{lemma}
Let $H < G$ be Nash groups and $M$ be a Nash $G$-manifold. Suppose
that the actions of $H$ on $G$ and of $G$ on $M$ are strictly
simple. Then the action of $H$ on $M$ is also strictly simple.
\end{lemma}
\emph{Proof.}\\
Consider the locally trivial fibration $M \rightarrow G \setminus
M$. If it is trivial, the statement is clear. It is locally
trivial and the statement is local.  \proofend \\
\emph{Proof of Proposition \ref{takoyto}.}

Case 1. $dimH=dimG$ \\
From the theory of Lie groups we know that in this case $H$ is a
union of connected components of $G$. $G$ has a finite number of
connected components by Proposition \ref{fincomp}. Hence $G/H$ is
finite.

Case 2. $H$ and $G$ are Zarisky closed in $GL_n$. \\
In this case they are linear algebraic groups, and for them this
statement is known.

Case 3. $G$ is Zarisky closed in $GL_n$.\\
Denote by $\overline{H}$ the Zarisky closure of $H$. It has the
same dimension as $H$ by Theorem \ref{dimZarclos}. From case 1 the
action of $H$ on $\overline{H}$ is strictly simple. From the case
2 the action of $\overline{H}$ on $G$ is strictly simple. Hence by
the lemma the action of $H$ on $G$ is strictly simple.

Case 4. General.\\
From the proof for case 1 we see that $G/H$ is a union of a finite
number of connected components of  $\overline{G}/H$ which is a
Nash manifold by case 3. \proofend
\section{Possible extensions and applications} \label{Summary}
We believe that it is possible to obtain an alternative proof of
de-Rham theorem which will be valid also in non-affine Nash case.
That proof goes in the following way. First one should prove for
$M = \R^n$ in the same way as we did. Then one should prove that
the cohomologies of a Nash manifold in classical topology are
equal to its cohomologies in the restricted topology and to
cohomologies of its de-Rham complex with generalized Schwartz
coefficients. If $M$ has a finite cover by open semi-algebraic
subsets Nash diffeomorphic to $\R^n$ such that all their
intersections are also Nash diffeomorphic to $\R^n$ then the
statement is easy because all these cohomologies are isomorphic to
the cohomologies of the \Che complex of this cover. But in general
the intersection of the open sets in the cover can be not Nash
diffeomorphic to $\R^n$. However we can always construct a
hypercover by open semi-algebraic sets Nash diffeomorphic to
$\R^n$. So one should prove that the \Che cohomologies of this
hypercover are isomorphic to the required cohomologies. For the
notion of hypercover see \cite{Del}.

After one proves de-Rham theorem for general Nash manifolds, the
relative de-Rham theorem and Shapiro lemma will follow in the same
way as in this paper.

It is possible to prove that for any Nash groups $H<G$, the action
of $H$ on $G$ is strictly simple. In fact, for any closed Nash
equivalence relation $R \subset M \times M$ we can build a
structure of $\R$-space on $M / R$. It is easy to see that if the
projection $pr:R \rightarrow M$ is \`{e}tale then $M / R$ is a
Nash manifold. It is left to prove that $M /R$ is Nash manifold in
case of any submersive $pr$. This problem is analogous to the
following known theorem in algebraic geometry. Let $M$ be an
algebraic variety. Let $R \subset M \times M$ be a closed
algebraic equivalence relation. Suppose that the projection $pr:R
\rightarrow M$ is smooth. Then $M / R$ is an algebraic space. This
theorem is proven using the fact that any surjective smooth map
has a section locally in \`{e}tale topology. In our case any
surjective submersion has a section locally in restricted
topology. So we think that our statement can be proven in the same
way.

In the classical case Shapiro lemma has a stronger version which
enables to compute cohomologies of $\mathfrak{g}$ in the case that
$G$ and $H$ are not cohomologically trivial. We think that our
techniques enable to prove its Schwartz version.

Using Shapiro lemma and \cite{AG} one can estimate
$H^i(\mathfrak{g}, \G(M,E))$, where $M$ is a Nash $G$ - manifold
with finite number of orbits, and $E$ is $G$-equivariant Nash
bundle over $M$. These cohomologies are important in
representation theory since sometimes the space of homomorphisms
between two induced representations is $H^0(G,\G(M,E))$ for
certain Nash bundle $E \to M$.

\appendix

\section{Proof of Theorem \ref{sursec}} \label{AppSurSec}
\setcounter{lemma}{0}

In this Appendix we prove Theorem \ref{sursec}. Let us first
remind its formulation.
\begin{theorem}
Let $M$ and $N$ be Nash manifolds and $\nu:M \rightarrow N$ be a
surjective submersive Nash map. Then locally it has a Nash
section, i.e. there exists a finite open cover $N= \bigcup \limits
_{i=1}^k U_i$ such that $\nu$ has a Nash section on each $U_i$.
\end{theorem}
This theorem follows immediately from the following three
statements.

\begin{theorem} \label{semisec}
Any semi-algebraic surjection $f:M \rightarrow N$ of
semi-algebraic sets has a semi-algebraic section.
\end{theorem}

\begin{theorem} \label{strat}
Let $f:M \rightarrow N$ be a semi-algebraic map of Nash manifolds.
Then there exists a finite stratification of $M$ by Nash manifolds
$M= \underset{i=1}{ \overset{k}{\udot}}  M_i$ such that $f|_{M_i}$
is Nash.
\end{theorem}

\begin{proposition} \label{ext}
Let $M$ and $N$ Nash manifolds and $\nu:M \rightarrow N$ be a Nash
submersion. Let $L \subset N$ be a Nash submanifold and $s:L
\rightarrow M$ be a section of $\nu$. Then there exist a finite
open Nash cover $L \subset \bigcup \limits _{i=1}^n U_i$ and
sections $s_i:U_i \rightarrow M$ of $\nu$ such that $s|_{L \cap
U_i} = s_i|_{L \cap U_i}$.
\end{proposition}

\subsection{Proof of Theorem \ref{semisec}}
\indent \indent Case 1. $M\subset N\times [0,1]$, $f$ is the
standard projection.
\\
We fix here a certain well-defined semi-algebraic way to choose a
section. One could do it in lots of different ways. For any $y\in
N$ define $F_y:=p(f^{-1}(y))$ where $p:M \rightarrow [0,1]$ is the
standard projection. $F_y\subset [0,1]$ is a semi-algebraic set,
hence a finite union of intervals. Let $\overline{F_y}$ be its
closure in the usual topology. Denote $s_1(y) := \min
\overline{F_y}$. Note that $s_1(y)$ is an end of some interval in
$F_y$. Denote this interval by $I_y$. Let $s_2(y)$ be the center
of $I_y$. Now define $s(y) := (y,s_2(y))$. By Seidenberg-Tarski
theorem $s$ is semi-algebraic, and it is obviously a section of
$f$.

Case 2. $M\subset N\times \R$, $f$ is the standard projection.
\\ We semi-algebraically embed $\R$ into $[0,1]$ using the
stereographic projection and reduce this case to the previous one.

Case 3. For $M\subset N\times \R^n$, $f$ is the standard
projection.
\\Follows by induction from case 2.

Case 4. General case. Follows from case 3 by considering the graph
of $f$. \proofend

\subsection{Proof of Theorem \ref{strat}}
In order to prove this theorem, we will need the following two
theorems from \cite{BCR}.

\begin{theorem}[Sard's theorem] \label{SSard}
Let $f:M \rightarrow N$ be a Nash map of Nash manifolds. Then the
set of its critical values is a semi-algebraic subset in $M$ of
codimension 1.
\end{theorem}
The proof is written on page 235 of \cite{BCR} (theorem 9.6.2).

\begin{theorem}[Nash stratification] \label{Nashstrat}
Let $M \subset \R^n$ be a semi-algebraic set. Then it has a finite
stratification by Nash manifolds $M = \udot N_i$.
\end{theorem}
The proof is written on page 212 of \cite{BCR} (theorem
9.1.8).\\\\
\emph{Proof of Theorem \ref{strat}}. It easily follows by
induction from the last two theorems and the following
observation. Let $f:M \rightarrow N$ be a semi-algebraic map
between Nash manifolds. Suppose that the graph $\Gamma_f$ of $f$
is a Nash manifold. Then the set of irregular points of $f$ is
exactly the set of critical values of the standard projection
$p:\Gamma_f \rightarrow M$. \proofend

\subsection{Proof of Proposition \ref{ext}}
\begin{notation}
{  Let $x\in \R^n$, $r\in \R$. We denote by $B(x,r)$ the open ball
with center $x$ and radius $r$. }
\end{notation}
\begin{definition}
{  A  Nash map $e:M \rightarrow N $ is called \textbf{\`{e}tale}
if for any $x \in M$, $de_x:T_xM \rightarrow T_{e(x)}N$ is an
isomorphism. }
\end{definition}
%
We will need a lemma from
\cite{AG} (Theorem 3.6.2).
\begin{lemma}
Let $N \subset \R^n$ be an affine Nash manifold and $L \subset N$
be a Nash submanifold. Then there exists a Nash positive function
$f_L^N:L \rightarrow \R$ and a Nash embedding $\phi_L^N:U_{f_L^N}
\hookrightarrow N$ such that $\phi(x,0) = x$, where $U_f :=
\{(x,y) \in N_L^N | ||y|| < f(x)\}$ and $||y||$ is the norm
induced from $\R^n$ to the normal space at $x$.
\end{lemma}
\emph{Proof of the proposition}.\\
Warning: proofs for cases 1 and 2 are technical and boring. The
reader will suffer less if he will do them himself.

Case 1. The map $\nu$ is \`{e}tale.\\
It is enough to prove for affine $M$ and $N$. Embed $M\subset
\R^k$ and $N \subset \R^l$. Consider the graphs $\Gamma(\nu)
\subset M \times N$ and $\Gamma(s) \subset \Gamma(\nu)$. Note that
$N_{\Gamma(s)}^{\Gamma(\nu)}$ is naturally embedded to
$\R^{2(k+l)}$ From differential topology we know that for any
$y\in M$ there exists $r\in \R$ such that $\nu|_{B(y,r)\cap M}$ is
an embedding. For any $((m,n),v) \in N_{\Gamma(s)}^{\Gamma(\nu)}$
denote $B_{((m,n),v)}(r):=B(((m,n),v),r)\cap
N_{\Gamma(s)}^{\Gamma(\nu)}$. Consider the function $g: \Gamma(s)
\rightarrow \R$ defined by $g(m,n) = \sup \{r\in \R |(pr \circ
\phi_{\Gamma(s)}^{\Gamma(\nu)}) |_{B_{((m,n),0)}(r)} $ is an
embedding $\}$/2,where $pr:\Gamma(\nu) \rightarrow N$ is the
standard projection. Denote
$h=min(\nu_{\Gamma(s)}^{\Gamma(\nu)},g)$. It is easy to see that
$\phi_{\Gamma(s)}^{\Gamma(\nu)}(U_h)$ is the graph of the required
section.

Case 2. $N \subset \R^l$ is affine,
 $M\subset \R^k \times N$ open, and $\nu$ is
the standard projection.\\
Consider the function $g: L \rightarrow \R$ defined by $g(x) =
\sup \{r\in \R | B(s(x),r) \cap N \times \R^k \subset M \}/2$. For
any $x\in L$ define $B_x = \nu(B(s(x),g(x))\cap M)$. For any
$(x,v) \in N_L^N$ define $B_{(x,v)}(r):=B((x,v),r)\cap N_L^N$.
Define $g_2:L \rightarrow \R$ by $g_2(x) = \sup \{r\in \R |
\phi_L^N(B_{(x,0)}(r)) \subset B_x\}/2$. Denote
$h=min(\nu_L^N,g_2)$. Now we define $s':\phi_L^N(U_h) \rightarrow
M$ by $s'(x)= (p(s(\pi((\phi_L^N)^{-1}(x)))),x)$, where $p:\R^k
\times N \rightarrow \R^k$ is the standard projection, $\pi:CN_L^N
\rightarrow L$ is the standard projection.

Case 3. For  $N \subset \R^l$ affine,
 $M\subset \R^k \times N$ any Nash submanifold, and $\nu$ is
the standard projection.\\
Denote $m:=dim(M)$ and $n:=dim(N)$. Let $\kappa$ be the set of all
coordinate subspaces of $\R^k$ of dimension $n-l$. For any $V \in
\kappa$ consider the projection $p:M \rightarrow N \times V$.
Define $$U_V:= \{x \in M | dp_x \text{ is an isomorphism }\}.$$ It
is easy to see that $p|_{U_V}$ is \`{e}tale and $\{U_V\}_{V \in
\kappa}$ gives a finite cover of $M$. Now this case follows from
the previous two ones.

Case 4. General case.
\\It is enough to prove for affine $M$ and $N$. Now we can replace $M$ by $\Gamma(\nu)$ and reduce to case 3.
\proofend

\end{document}